\documentclass[secnumarabic,aps,showpacs,showkey,superscriptaddress,floatfix]{revtex4-1}

\usepackage{amsmath,amssymb,amsthm}

\usepackage{hyperref}
\usepackage{amsthm}
\usepackage{graphicx}
\usepackage{subfigure}
\usepackage{color,bbm}
\usepackage[normalem]{ulem}

\definecolor{violet}{rgb}{1.00,0.00,1.00}	

\graphicspath{{./Figures/}}

\newcommand{\R}{\mathbb{R}}

\newcommand{\C}{\mathbb{C}}

\renewcommand{\P}{\mathcal{P}}

\begin{document}

\title{Index Distribution of the Ginibre Ensemble }

\author{Romain Allez}%
\email[]{romain.allez@gmail.com}
\affiliation{Weierstrass Institute, Mohrenstr. 39 · 10117 Berlin, Germany. }
\author{Jonathan Touboul}
\email[]{jonathan.touboul@college-de-france.fr}
\affiliation{The Mathematical Neuroscience Laboratory, CIRB / Coll\`ege de France (CNRS UMR 7241, INSERM U1050, UPMC ED 158, MEMOLIFE PSL*)}
\affiliation{SISYPHE Laboratory, INRIA Paris}
\author{Gilles Wainrib}
\email[]{wainrib@math.univ-paris13.fr}
\affiliation{LAGA, Universit\'e Paris XIII}

\date{\today}%

\begin{abstract}
Complex systems, and in particular random neural networks, are often described by randomly interacting dynamical systems with no specific symmetry. In that context, characterizing the number of relevant directions necessitates fine estimates on the Ginibre ensemble. In this Letter, we compute analytically the probability distribution of the number of eigenvalues $N_R$ with modulus greater than $R$ (the index) of a large 
$N\times N$ random matrix in the real or complex Ginibre ensemble. We show that the fraction $N_R/N=p$ has a distribution scaling as $\exp(-\beta N^2 \psi_R(p))$ with $\beta=2$ (respectively $\beta=1$) for the complex (resp. real) Ginibre ensemble. For any $p\in[0,1]$, the equilibrium spectral densities as well as the rate function $\psi_R(p)$ are explicitly derived. This function displays a third order phase transition at the critical (minimum) value $p^*_R=1-R^2$,  associated to a phase transition of the Coulomb gas. We deduce that, in the central regime, 
the fluctuations of the index $N_R$ around its typical value $p^*_R N$ scale as $N^{1/3}$.    
\end{abstract}

\maketitle

In such different domains as physics, biology, economics and social science, large disordered networks describing randomly interacting dynamical agents are ubiquitous~\cite{boccaletti2006complex}. Randomness in the interactions have made of Random Matrix Theory (RMT) a central tool in the analysis of these systems. An emerging application of RMT is the theory of random neuronal networks of the brain. In contrast with most physical systems described as random Hamiltonian (such as spin glasses) where the interconnections matrices are symmetrical, cortical interconnections are asymmetric and display specific disorder levels~\cite{parker:03,marder-goaillard:06} resulting from development and learning. A generic example used both in theoretical neuroscience and machine learning~\cite{amari:72,jaeger:04} is given by the neural network model describing the activity $\mathbf{x}$ of $N$ neurons through the equations:
\begin{equation}\label{eq:SompoDiscrete}
\mathbf{x}(t+1)= S(\sigma \mathbf{J}\cdot\mathbf{x}(t))
\end{equation}
for $\mathbf{J}\in \R^{N\times N}$ a random matrix and $S$ a centered sigmoidal function.

Similarly to disordered physical systems such as spin glasses or liquids which display a rich (free) energy landscape characterized by many extrema and complex stability patterns~\cite{wales}, random neural networks display explosion of complexity as a function of the disorder strength $\sigma$~\cite{wainrib:13}. 
Similar models for stability of ecosystems were studied by May~\cite{may} in the case of a symmetric interconnection matrix $\mathbf{J}$: a sharp transition of the probability of stability in the system (weakly/strongly interacting phase) 
is shown in the limit of infinite dimensions. This sharp transition was then further investigated 
for large but finite $N$, which led to many developments in random matrix theory, in particular on the 
statistics of the maximal eigenvalues of random symmetric matrices (see \cite{review} for a review on this subject). 

At this transition, the behavior of the system \eqref{eq:SompoDiscrete} in a neighborhood of an equilibrium is characterized by its index, namely the number of unstable directions, that are the relevant directions of the dynamics~\cite{lajoie2013chaos}. The purpose of this paper is specifically to characterize the index around fixed points of discrete-time random dynamical systems. This quantity is given by the number of eigenvalues of the Jacobian (non-symmetric) random matrix that have a modulus larger than $1$. In the model 
\eqref{eq:SompoDiscrete}, the trivial equilibrium has therefore an index equal to the number of eigenvalues of $J$ with modulus larger than $1/\sigma$. Thus, a natural question one faces is: what is the probability of having $N_{R}$ eigenvalues outside the open ball 
$\mathcal{B}(0,R)=\{z\in \C: |z|<R \}$? 

A similar question can be raised in the analysis of Hamiltonian systems. In contrast with our case, the Hessian matrix of the Hamiltonian characterizes the stability of an equilibrium. This matrix is Hermitian, hence has real eigenvalues. A one-dimensional Coulomb gas approach was introduced for Gaussian and Wishart ensembles~\cite{majumdar-nadal:09,majumdar-nadal:11,majumdar-vivo:12} to address these questions. 

Much less is know about the real or complex Ginibre ensembles~\cite{ginibre1965statistical} consisting of random matrices with independent Gaussian entries with variance $1/{N}$. Of course, such matrices are not symmetric and non-normal (i.e., do not commute with their Hermitian conjugate).
They are nevertheless diagonalizable in $\C$ and 
the joint probability density (jpdf) function of the eigenvalues~\cite{bordenave-chafai:12} is, in the complex case,  
\begin{equation}\label{jpdf}
P(z_1,\cdots,z_N) = \frac{1}{Z_N} \exp\left(- N \sum_{k=1}^N |z_k|^2 + 2 \sum_{i< j} \ln |z_i-z_j|\right) 
\end{equation}
where $Z_N$ is a normalization constant. The joint distribution of the real Ginibre ensemble is slightly different due to the symmetry of the spectrum with respect to the complex conjugation~\cite{lehmann1991eigenvalue} and does not have a density 
with respect to Lebesgue in $\C^N$.
In physical terms, the jpdf \eqref{jpdf} of the eigenvalues of the complex Ginibre Ensemble 
describes the equilibrium Boltzmann weight of 
 a 2-dimensional Coulomb gas of $N$ positively charged particles with electrostatic repulsion, 
 confined in a quadratic potential and subject to a thermal noise~\cite{forrester}. 
 
\begin{figure}
\centering
\includegraphics[width=.65\textwidth]{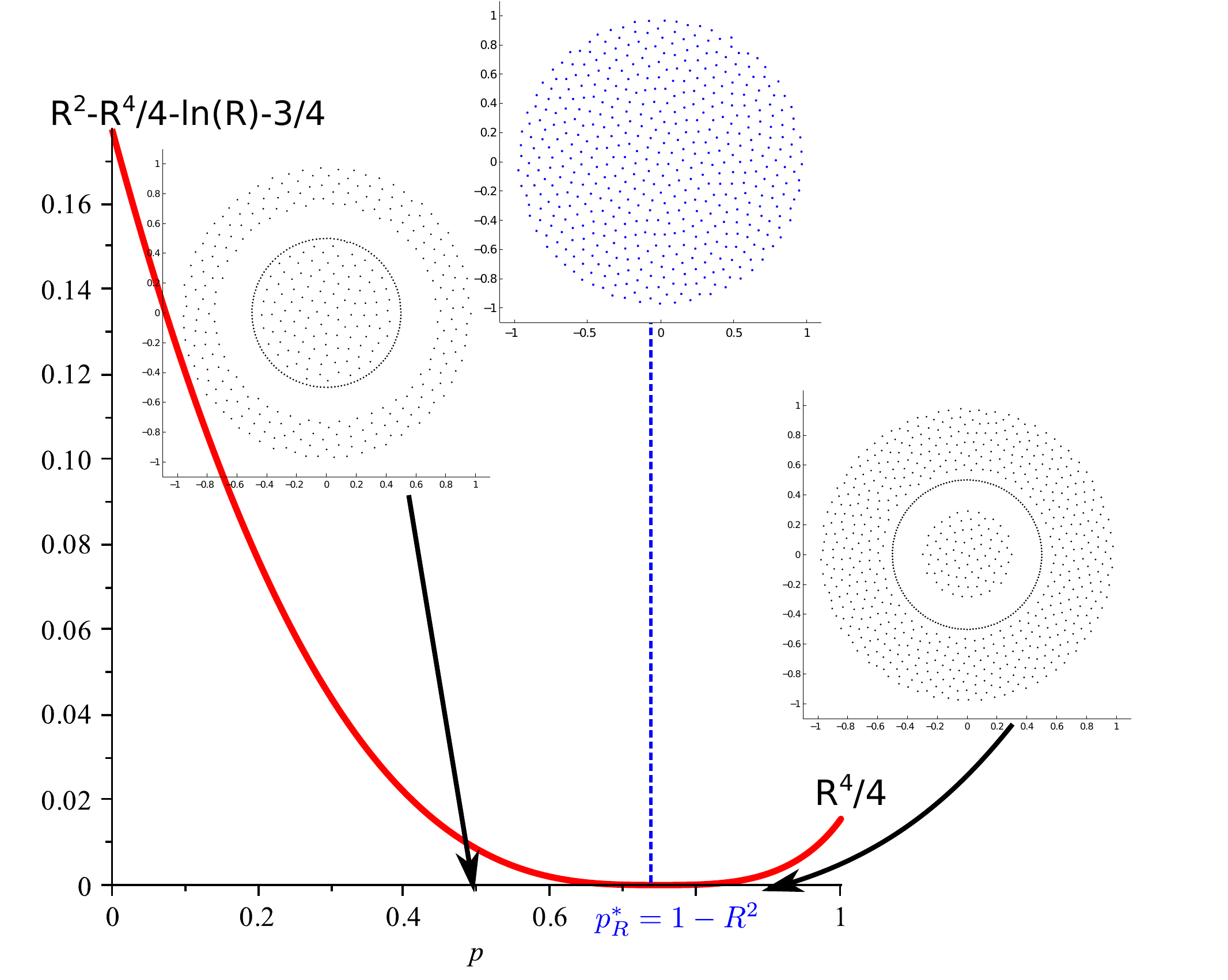}
\caption{Rate function $\psi_R(p)$ for $R=1/2$. The function shows a minimum at $p=p^*(R)=1-R^2$ where the phase 
	transition for the empirical density of charges occurs. The insets display scatter plots of numerically 
	simulated Coulomb gas under the constraint $\{N_R=p N\}$
	for $p=1/2 < p^*(R)$ (left), $p=p^*(R)$ (middle) and $p=9/10>p^*(R)$ (right). }
\label{fig:Psi}
\end{figure}

Using large deviations techniques (saddle point method), we compute analytically, for any $p\in[0,1]$, the probability $\P(p,N)$ that a fraction $p=N_R/N$ of eigenvalues lie outside $\mathcal{B}(0,R)$ in the large $N$ limit.
The method relies on the large deviation principles satisfied by the empirical spectral
density for both the complex and real Ginibre ensembles. 
In the complex case, the eigenvalues 
$z_1,z_2,\dots,z_N$ of a complex $N \times N$ Ginibre matrix are distributed according to \eqref{jpdf} and their empirical measure is 
defined as
\begin{align*}
\rho_N(z) = \frac{1}{N} \sum_{k=1}^N \delta(z-z_k)\,. 
\end{align*}
Although the eigenvalues of a real $N \times N$ Ginibre matrix do not carry the same number of degrees of freedom as the corresponding complex matrix, the empirical measure in the real case satisfies the same large deviation principle as in the complex case.
Using those large deviation principles together with geometrical considerations on the distribution of particles in the Coulomb gas, we derive analytically the asymptotic spectral density of Ginibre ensembles constrained to have an asymptotic index $N_R\sim p N$. This enables to derive analytically the rate function $\psi_R(p)$ such that, at leading order:
\begin{align}\label{eq.general}
\P(p,N) \approx \exp\left[- \frac{\beta}{2} \, N^2 \, \psi_R(p) \right]\,, \quad 0 \le p \le 1
\end{align}
with $\beta=2$ (resp. $\beta=1$) in the complex (resp. real) case and
where the symbol $\approx$ denotes the logarithmic equivalent. The rate function $\psi_R(p)$ has a particularly simple form~\eqref{eq:Psi}. It is convex (see Fig.~\ref{fig:Psi}), and as expected, reaches its minimum (zero) at $p^*_R=1-R^2$ corresponding to the proportion of eigenvalues outside $\mathcal{B}(0,R)$ in the circular law. In contrast with symmetric cases~\cite{majumdar-nadal:09,majumdar-vivo:12}, the second derivative of $\psi_R$ continuously vanishes at the 
critical value $p^*_R$. The typical Gaussian regime of the fluctuations valid for one-dimensional Coulomb gases is therefore lost, as well as the typical size of the fluctuations that were found of order $\sqrt{\ln(N)}$. Here, the second derivative has instead an {\it absolute value singularity}, i.e. the third derivative $\psi_R^{(3)}$ is discontinuous at $p^*_R$:
\begin{align}\label{cubic.behavior}
\psi_R(p^*_R + \delta) = \frac{|\delta|^3}{6R^2} +o (\delta^3). 
\end{align} 
This cubic behavior implies that the central fluctuations of the index $N_R$ around its typical value $p^*_R N$ for large $N$ are much wider than in one-dimensional gases: the index $N_R$ fluctuates in the central regime around its typical value $(1-R^2) N$ on a region of scale $N^{1/3}$. Interestingly, this {\it third order} phase transition is associated to a phase transition of the Coulomb gas at $p=p^*_R$, as illustrated in Fig.~\ref{fig:Psi} where are plotted the empirical eigenvalue densities found in the two regimes $p<p^*_R$ and $p>p^*_R$, computed here through direct simulations of the constrained Coulomb gas.

Let us now go into the details of the derivation. Both in the complex~\cite{petz} and the real~\cite{gerard} Ginibre ensembles, the empirical eigenvalue density $\rho_N$ satisfies a {\it large deviation principle} with speed $N^2$ and rate function
\begin{align}\label{rate.function}
\mathcal{I}[\mu] := \frac{\beta}{2} \left (\int_{\C} |z|^2 \, \mu(dz) - \int_{\C^2} \ln|z-z'| \mu(dz) \mu(dz')-\frac 3 4 \right)\,. 
\end{align} 
where $\beta=2$ (resp. $1$) for the complex (resp. real) case.  
This means that if $\mathcal{A}$ is a subset of the set $\mathcal{M}_1(\C)$ of probability measures on $\C$, then 
\begin{align}\label{largedev}
\P[\rho_N\in \mathcal{A}] \approx \exp\left[- N^2 \,\min_{\mathcal{A}} \mathcal{I} \right]\,. 
\end{align} 
The derivation of the asymptotic eigenvalue density under some constraint therefore consists in minimizing the rate function $\mathcal{I}[\mu]$ for $\mu$ satisfying the constraint. Once this minimizer $\mu^*$ is known, the probability of the event is deduced from Eq. \eqref{largedev}  by computing the value of the rate function at its minimizer.

Characterizing the index large deviations thus amounts to finding the equilibrium density of the Coulomb gas constrained to have $N_R\sim pN$ particles outside $\mathcal{B}(0,R)$. This is generally a technically involved task. For one-dimensional gases for instance~\cite{majumdar-nadal:09,majumdar-nadal:11,majumdar-vivo:12,majumdar-vergassola:09}, this optimization problem necessitates to solve a singular integral equation involving Cauchy's principal part. Here the problem is strikingly simplified by the 2-d nature of the associated gas as well as the rotational invariance of the system, which impose a specific and simple form for the possible minimizers. 

The limiting distribution $\mu^*$ indeed preserves the invariance under rotation satisfied by both the rate function~\eqref{rate.function} and the constraint considered here. It can therefore be decomposed in a sum of radial probability densities supported on rotation invariant subsets of $\C$. Such sets are either balls, annuli or circles centered at the origin. The invariance under rotation imposes $\mu^*$ to have a uniform density along the circles. 
Furthermore, it is easy to show that the measure $\mu^*$ has a constant density 
equal to $1/\pi$ on the interior of its support. 
Indeed, differentiating $\mathcal{I}[\mu] $ at $\mu^*$ yields for any $z\in \C$
\begin{equation}\label{eq.minimizer}
|z|^2 - 2 \int_{\mathbb S} \ln|z-z'|\, \mu^*(dz') = 0 
\end{equation}
where $\mathbb S $ is the support of $\mu^*$. 
Taking the Laplacian with respect to $z$ in the interior of the support of $\mu^*$ 
and using the fact that the logarithm is the Green function of the Poisson equation on the plane, 
we obtain that $\mu^*(z)=1/\pi$. We recover here a general result~\cite{asz:13}: the limiting empirical eigenvalue density decomposes into a sum of regular measures with a constant density $1/\pi$ and singular measures supported on the boundary of the constraint (respective supports, and therefore shape of the minimizers, are generally unknown). 

Of course, under no constraint, the distribution $\mu^*$ is regular on its full support, and we recover the {\it circular law}:
$\mu^*(dz)= \frac{1}{\pi}1_{\{|z| \le 1\}} \, dz$
which corresponds  to $p=p^*_R$ (the radius of the ball is found as the unique value ensuring normalization of $\mu^*$). 

For $p\neq p^*_R$, the constraints are active and necessarily, the measure has a singular part supported by centered circles. The minimization of the rate function is now reduced to finding the radii of the ball, annuli and circles supporting the density, which can be done analytically in all cases. Let us treat for example the case where we impose an atypically small proportion of charges $p<1-R^2$ to 
lie outside the open ball ${\mathcal{B}}(0,R)$ (i.e. an anomalously large proportion of charges $1-p>R^2$ in the closed ball 
$\overline{\mathcal{B}}(0,R)$).
As observed, the density of the constrained equilibrium measure is $1/\pi$ inside the 
open ball $\mathcal{B}(0,R)$, which accounts for a proportion equal to $R^2$ of the charges. The excess proportion $(1-p)-R^2$ of charges which is pulled inside $\mathcal{B}(0,R)$ by the constraint have no choices but lying exactly on the circle of radius $R$, on which they allocate uniformly.  Outside $\overline{\mathcal{B}}(0,R)$, we expect the charges to allocate in one single annulus~\footnote{Note that here, we ruled out the possibility of having a ball of radius strictly smaller than $R$ or several disjoint annuli. Such measures are not physically sound. Analytically, our approach allows to compute and optimize on the radius of the ball and circle and find that both radii are equal to $R$, and check that if one considers two annuli, these will necessarily aggregate into one annulus.}. 
Functions achieving the minimum of the rate function are therefore of the form:
\begin{align*}
\mu^*(dz) = \frac{1}{\pi} 1_{\mathcal{B}(0,R)}(z) \, dz + \frac{1-p-R^2}{2\pi}\, \delta(r-R) \, d\theta\, + 
\frac{1}{\pi} 1_{\mathcal{A}(R_1,R_2)}(z) \, dz,
\end{align*}
where $\mathcal{A}(R_1,R_2)=\{z\in\C; R_1< \vert z \vert < R_2\}$ and $R_2^2-R_1^2 =p$ by normalization of $\mu^*$. 
The problem of minimizing $\mathcal{I}[\mu^*]$ reduces to a single variable optimization on $R_1\in [R,\infty)$. This is performed analytically by computing explicitly $\mathcal{I}[\mu^*]$ using the identity 
valid for any $z,z'\in \C\setminus \{0\}$
\begin{align}\label{identity}
\int_0^{2\pi} \ln |z- z' e^{i\theta}| \, d\theta = 2\pi \, \ln(\max(|z|, |z'|))
\end{align}
which arises from the harmonic nature of $\ln |z|$ on $\C\setminus \{0\}$. We eventually obtain 
\begin{align*}
\mathcal{I}(\mu^*) = \frac{\beta}{2}\bigg[- \frac{R^4}{4} + &(1-p) R^2 - (1-p) ^2 \ln R  + \frac{3}{4} p^2 +
 \frac{p}{2}  R_ 1^2 + p (1-p)
- R_1^2(R_1^2-2(1-p) ) \ln R_1  \\ 
&+\frac{1}{2} \left((R_1^2-1+p )^2 -1 \right) \ln(p+R_1^2) \bigg] =: J_R(R_1)\,. 
\end{align*}
It is now straightforward to compute the derivative $J_R'(R_1)= \beta R_1(R_1^2-1+p) \ln \frac{p+R_1^2}{R_1^2}$
and to check that $J_R$ reaches its minimum on the interval $[R,+\infty)$ for $R_1=\sqrt{1-p}$. 
Finally the equilibrium measure under the constraint is 
\begin{align}\label{minimizer}
\mu^*(dz) = \frac{1}{\pi} 1_{\mathcal{B}(0,R)}(z) \, dz + \frac{1-p-R^2}{2\pi}\, \delta(r-R) \, d\theta\, + 
\frac{1}{\pi} 1_{\mathcal{A}(\sqrt{1-p},1)}(z) \, dz\,.
\end{align}
This is precisely the measure found numerically by simulating the constrained Coulomb gas system
(see the left scatter plot of Fig. \ref{fig:Psi}). The case $p>1-R^2$ is handled similarly and leads to the probability measure with a ball of radius $\sqrt{1-p}$, an annulus of inner and outer radii $(R,1)$ with uniform distribution $1/\pi$, the rest of the mass being uniformly distributed along the circle of radius $R$. This result is again in agreement with our numerical simulations (see the right inset of Fig. \ref{fig:Psi}). The two different shapes of the density prove that the Coulomb gas undergoes a phase transition at $p=p^*_R$. 

Finally, the large deviation estimates of the index $N_R$ are found simply by evaluating 
the rate function $\mathcal{I}$ at its minimizer $\mu^*$; we get  
\begin{equation*}
\P(p,N) \approx \exp\left(- \frac{\beta}{2} \, N^2 \, \psi_R(p)\right)
\end{equation*}
where $q=1-p$ and 
\begin{equation}\label{eq:Psi}
	\psi_R(p)=\begin{cases}
		\frac{R^4}{4} - q R^2 + q^2 (\frac{3}{4}-\frac{\ln q}{2} + \ln R)\,, & p>1-R^2\\
		 -\frac{R^4}{4} + q R^2 -q^2 (\frac{3}{4}   - \frac{\ln q}{2} + \ln R)\,, & p<1-R^2\,. 
	\end{cases}
\end{equation}
It is easy to see from this simple formula the discontinuity of the third derivative announced~\eqref{cubic.behavior}. The probability of the event $\{N_R\sim (p^*_R+\delta)N \}$ is non-degenerate for $\delta$ of order $N^{-2/3}$, meaning that the index $N_R$
fluctuates in the central regime around its typical value $(1-R^2) N$ on a region of width $N^{1/3}$. This scale contrasts with the case of symmetric Gaussian $\beta$ and Wishart ensembles where the typical fluctuations spread on a much narrower region of order $\sqrt{\ln N}$. The dissimilar extents of the central regimes is evocative of the fact that electrostatic 
repulsion has stronger effects in one dimensional than in two dimensional gases.

\begin{figure}
	\centering
		\includegraphics[width=.75\textwidth]{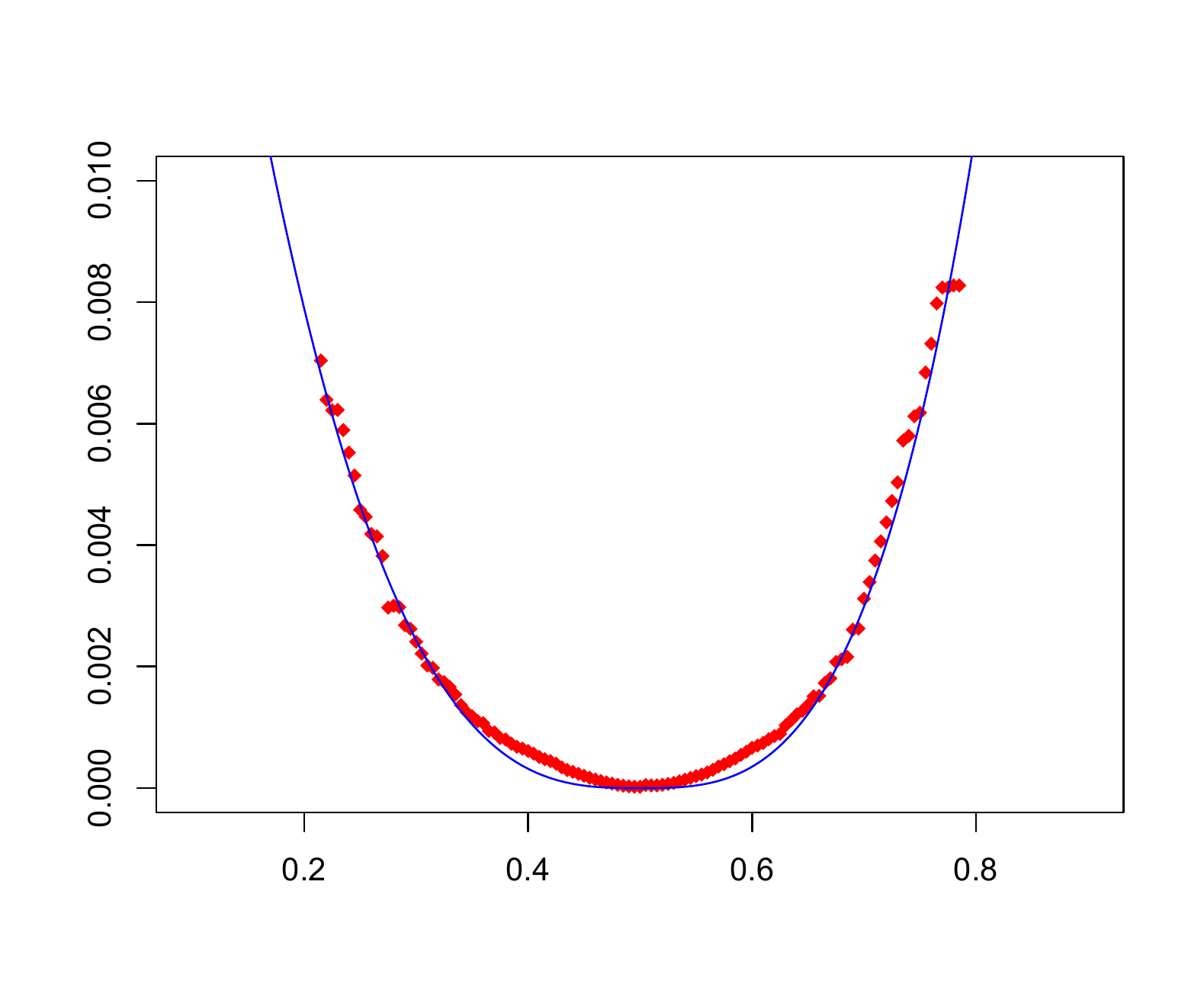}
	\caption{(Color online). The blue curve represents the theoretical 
	rate function $\psi_R(p)$ \eqref{eq:Psi} in the symmetric regime ($R=1/\sqrt{2}$). The red dots display $-\ln \P(p,N)/N^2$
	for $p=i/N, i=1,\dots,N$ computed from our sample constructed using modified Metropolis algorithm with $10^5$ iterations and   $N=200$ particles.}
	\label{fig:FitsReal}
\end{figure}

We confirm the analytical prediction \eqref{eq:Psi} numerically. In order to explore the extremely rare events characterized by $\psi_R$, 
a particularly efficient method is the Metropolis algorithm in its modified version as presented in~\cite{majumdar-nadal:11} taking into account conditioning and the discrete nature of possible values of $N_R$. 
The obtained numerical curve fits our theoretical large deviation rate function as given in \eqref{rate.function} 
with very good agreement (see Fig. \ref{fig:FitsReal}) away from the central regime. 

Characterizing the fluctuations of the index in the central regime remains an exciting open problem. Large deviations techniques 
as used here do not permit to understand those typical fluctuations. Nevertheless
we expect the central regime to establish a smooth crossover between the left and right large deviation regimes 
as it was found for the top eigenvalue of Gaussian Wigner and Wishart matrices in \cite{majumdar-vergassola:09} 
where the authors exhibit matching of the Tracy Widom tails with the two left and right large deviation rate functions near 
the scaling region. 
A correct matching in our case would force  
the central regime probability distribution to 
have symmetric tails decaying at leading order as $\exp(-|x|^3/(6R^2))$ for $x\to\pm \infty$. 

Our results \eqref{minimizer} and \eqref{eq:Psi} 
can be extended straightforwardly in the case of general $d$-dimensional Coulomb gas with $d\geq 2$.
The equilibrium jpdf of the $d$-dimensional Coulomb gas expresses in term 
of the Green function $G_d$ in dimension $d$ as 
\begin{align*}
P(z_1,\dots,z_N) = \frac{1}{Z_N} \exp\left(- N  \sum_{i=1}^N |z_i|^2 +2\sum_{i< j} G_d(z_i-z_j)\right)\,. 
\end{align*} 
Note that the repulsion between the particles is no longer logarithmic if $d>2$ since in this case, $G_d(z)= 1/|z|^{d-2}$. 
The two dimensional case remains somehow the most interesting to our eyes and for the applications  
because of its direct relation with random matrices and applications to theoretical neurosciences. 

Questions related to the present study were addressed in \cite{lebowitz}
where the authors investigate macroscopic large deviations in the two and three dimensional plasmas.
The methods used are very different from ours, and the results of the present manuscript recover a limit case treated in~\cite{lebowitz}, Eq. $(2.7)$, corresponding to our formula \eqref{eq:Psi} in the particular case $p=1$ and $R>0$.   

Back to our study, a further interesting question is to compute the mean 
empirical density of charges $\langle \rho_N(z)\rangle$ for fixed values of $N$ and under the above constraints. 
For any large but finite $N$, $\langle \rho_N(z)\rangle$  is smooth with respect to $z$ and closely approximates 
the singular probability measure $\mu^*(dz)$ which appears at $N=+\infty$. 
We expect this finite $N$ deformation, and in particular its scaling shape
near the circle which supports the limiting singular density, to have strong implications for the dynamic of random neural networks (see below). 
In this direction, a tentative method is to account for the correction $\rho^*_1(z)$ 
such that $\langle\rho_N(z)\rangle= \mu^*(dz)+\rho^*_1(z)/N +o(1/N)$, 
induced by the entropic term of order $1/N$ in the rate function \eqref{rate.function}, 
(this correction was computed explicitly for the Gaussian \cite{jp-alice} and Wishart-Laguerre \cite{satya} symmetric ensembles 
under no constraint). 
The main difficulty here is that the entropic term makes 
the optimization of the rate function \eqref{rate.function} much harder to analyze as 
one can no longer postulate a specific form for the possible
minimizers. Moreover this {\it perturbative} method~\footnote{Under no constraint (free case), this study has been worked 
out in \cite{preparation} but remains open under some constraint.} is only 
valid {inside} the support and breaks down near the singularities.

Strikingly, we note that the constraint has a particularly simple effect on the limiting density of charges in the Coulomb gas, evocative of an optimal transport of mass problem. Indeed, the minimal ``effort'' to take the circular distribution to fulfill the constraint simply amounts to displace the closest charges onto the circle of radius $R$. The minimizer of the constrained problem corresponds to the \emph{projection} of the circular law on the set of measures with mass $p$ outside $\mathcal{B}(0,R)$ in the sense that it minimizes any transport cost compatible with the Euclidian distance~\cite{villani2009optimal}.
This reveals a nice coincidence between constrained Ginibre ensemble, optimal transport and projection of measures. Such a relation cannot arise in one-dimensional gases with logarithmic repulsion between the particles, 
because the singular interaction potential is not compatible with the presence 
of Dirac masses. Note however that in dimension one, the logarithmic 
interaction is {\it not} the Coulomb interaction and this is certainly the reason why there is no 
``condensation" (Dirac mass on a singular point) phenomenon in the one dimensional case (see again 
\cite{majumdar-nadal:09,majumdar-nadal:11,majumdar-vivo:12}). If one considers the true 1D-Coulomb interaction,
$\propto |z_i-z_j|$ and not $\propto \ln |z_i-z_j|$, then one would also expect such condensation phenomenon to occur.

The results of this Letter also open exciting perspectives on the dynamics of random neural networks. 
Indeed, as disorder is increased in the connectivity matrix of the random neural network~\eqref{eq:SompoDiscrete}, the properties of the equilibria will undergo bifurcations that can be inferred from the study of Jacobian matrices. The local dynamics arising in the events considered here will modify the dynamics: at phase transition, decay rates and fluctuations will be highly sensitive to the spectrum of the coupling Ginibre matrix. Moreover, characterizing the dynamics necessitates to reduce the system on the center manifold, whose dimension is equal to the number of neutral directions of the dynamics (neither contracting nor expanding). Here, we have seen that prescribing an anomalous index does not affect the number of stable or unstable directions, but the excess or default of eigenvalues accumulate on the circle corresponding to neutral directions. In dynamical systems terms, such a conditioning increases the dimension of the center manifold, allowing very complex transitions to occur and rich dynamics.

{\bf Acknowledgments.} RA and JT warmly thank the organizers of the Beg Rohu summer school where the seeds of this work were sowed. RA is grateful to J.-P. Bouchaud, S.N. Majumdar and P. Vivo 
for enlightening discussions on saddle point analysis applied to Ginibre ensemble. JT and GW acknowledges K. Pakdaman for very interesting discussions on RMT and Coulomb gases. RA received funding from the European Research Council under the European
Union's Seventh Framework Programme (FP7/2007-2013) / ERC grant agreement nr. 258237 and thanks the Statslab in DPMMS, Cambridge for its hospitality at the time this work was finished.

\bibliographystyle{apsrev4-1} 
\bibliography{ginibre}

\end{document}